\documentclass[12pt]{article}
\usepackage[utf8]{inputenc}
\usepackage[T1]{fontenc}
\usepackage{amsmath,amssymb,amsthm}
\usepackage{geometry}
\usepackage{hyperref}
\usepackage{indentfirst}
\usepackage{booktabs}
\usepackage{longtable}
\usepackage{xcolor}
\usepackage{fvextra}
\allowdisplaybreaks

\DefineVerbatimEnvironment{code}{Verbatim}{
  breaklines=true,
  breakanywhere=false,
  fontsize=\fontsize{10}{10}\selectfont,
  xleftmargin=1em,
}
\DefineVerbatimEnvironment{term}{Verbatim}{
  breaklines=true,
  breakanywhere=true,
  fontsize=\fontsize{10}{10}\selectfont,
  xleftmargin=1em,
}

\newtheorem{theorem}{Theorem}
\newtheorem{lemma}{Lemma}
\newtheorem{proposition}{Proposition}
\newtheorem{corollary}{Corollary}
\theoremstyle{remark}
\newtheorem{remark}{Remark}
\theoremstyle{definition}
\newtheorem{definition}{Definition}

\title{Excluding 54 of the 76 Nontrivial Factorization Patterns of the Cuboid Polynomial $P_{a,b,u}(t)$}

\author{Valery Asiryan\\[3pt]
\small \texttt{asiryanvalery@gmail.com}}
\date{\small November 4, 2026}

\begin{document}
\maketitle

\begin{abstract}
Let $P_{a,b,u}(t)$ be Sharipov's even monic degree-$12$ third cuboid
polynomial, with positive pairwise distinct rational parameters.
We prove a uniform normal form for its complete factorization over
$\mathbb{Q}$. Every factorization consists of at most two even
quadratic factors, even irreducible factors of higher degree, and pairs
of irreducible factors exchanged by $t\mapsto-t$.
After dividing all degrees by the contribution of this symmetry, the
normal form satisfies one degree equation, one sign equation, and an
obligatory-factor alternative: an even irreducible factor of degree at
least six, or an even irreducible quartic with positive rational-square
constant term.

The arithmetic input is a uniform construction over $\mathbb{Q}_2$,
using Newton polygons, Hensel lifting, and a norm obstruction from
$\mathbb{Q}_2(i)$. Valuations at all finite primes constrain rational
factor constants to signed squares, while the real root equation gives
exactly two simple negative roots of the auxiliary sextic.
The normal form leaves exactly $22$ of the $76$ nontrivial partitions
of $12$, thereby excluding $54$ complete irreducible degree patterns
by one argument. Of these exclusions, $37$ follow from evenness alone
and $17$ from the additional arithmetic restrictions.
All local cases are proved explicitly, and an exact verification
script accompanies the algebraic identities and exhaustive partition
enumeration. The remaining patterns are necessary possibilities;
their realization is not asserted, and full irreducibility remains open.


\end{abstract}

\section{Introduction}\label{sec:intro}

The perfect cuboid problem asks for a rectangular box with positive
integer edges for which all three face diagonals and the space
diagonal are integers. Sharipov's polynomial approach leads to three
irreducibility conjectures~\cite{SharipovCuboids}. The third concerns
a degree-$12$ even polynomial depending on three positive integer
parameters; its exceptional parameter relations and its connection
with a Diophantine equation are described in~\cite{SharipovThird}.

\medskip
\noindent\textbf{Goal.}
We establish a common structural obstruction to complete
factorizations of the third polynomial. Instead of excluding degree
patterns individually, we obtain a normal form governed by the
involution $t\mapsto-t$, square classes of constant terms, and one
uniform local factor. The resulting conditions exclude $54$ of the
$76$ reducible degree partitions of twelve.

Evenness alone excludes $37$ of these partitions, since every odd
degree must occur an even number of times. The arithmetic restrictions
proved here reduce the remaining $39$ patterns to $22$.
Sharipov's structural theorem~\cite[Theorem~6.1]{SharipovThird}
characterizes integer roots through $48$ positive integer variables.
Our result instead constrains the complete irreducible factorization;
it does not exclude rational roots.

\medskip
\noindent\textbf{Method.}
The polynomial is a composition $P(t)=F(t^2)$, where $F$ is a sextic
given by a signed sum of three squared cubic expressions.
Its root equation implies that every monic rational factor of $F$
has constant term in $\mathbb{Q}^{\times2}\cup-\mathbb{Q}^{\times2}$.
Over the reals, the same equation has precisely two negative roots.
Over $\mathbb{Q}_2$, an exhaustive valuation argument supplies an even
irreducible factor of degree six or eight, or an even irreducible
quartic with square constant term. An even local factor forces its
containing rational irreducible factor to be even.
If that rational factor is quartic, its signed-square constant must
be positive. These inputs give a signed partition of six, from which
the complete degree inventory follows.

The arguments apply to every positive pairwise distinct rational
triple. In particular they apply in the ordered domain $0<u<a<b$.
They do not require coprimality or the additional exclusions
$a^2\ne bu$ and $b^2\ne au$. Thus the conclusions also constrain
the factorization on those two exceptional relations, without
claiming irreducibility there. No result on the first or second
cuboid conjecture is used.

Section~\ref{sec:poly} fixes the polynomial and states the normal-form
theorem. Sections~\ref{sec:global} and~\ref{sec:real} prove its
global and real inputs. Sections~\ref{sec:local}--\ref{sec:units}
give the full local construction. Section~\ref{sec:normalproof}
proves the normal form and its exclusion consequences.
Section~\ref{sec:excluded} then collects all $54$ excluded patterns in
Table~\ref{tab:excluded}, grouped by the condition that rules them out.
Only after this proved status has been recorded does
Section~\ref{sec:inventory} present the $22$ patterns not excluded by
the argument, in Table~\ref{tab:remaining}, and verify their formal
compatibility with the normal form. The appendix records a reproducible
exact verification script and its computed transcript.

\section{The polynomial and the normal-form theorem}\label{sec:poly}

Throughout the global argument, let $a,b,u\in\mathbb{Q}_{>0}$ be
pairwise distinct. Put
\begin{equation}\label{eq:ABU}
A=a^2,\qquad B=b^2,\qquad U=u^2.
\end{equation}
For arbitrary nonzero $A,B,U$, define
\begin{equation}\label{eq:F}
\begin{aligned}
F_{A,B,U}(x)={}&((x-A)(x+B)(x+U))^2\\
&+((x-B)(x+A)(x+U))^2\\
&-((x-U)(x+A)(x+B))^2.
\end{aligned}
\end{equation}
We abbreviate this polynomial to $F$ when the parameters are fixed.
It is monic of degree six. With $S=A+B$ and $T=AB$, its expansion is
\begin{align}
F(x)={}&x^6+(6U-2S)x^5
 +(S^2+4US+U^2-14T)x^4\nonumber\\
&+(6US^2-20UT-2S(U^2+T))x^3\nonumber\\
&+(T^2+U^2S^2-14TU^2+4UST)x^2\nonumber\\
&+(6UT^2-2U^2ST)x+U^2T^2.\label{eq:Fexpanded}
\end{align}
Sharipov's third polynomial is
\begin{equation}\label{eq:P}
P_{a,b,u}(t)=F_{a^2,b^2,u^2}(t^2).
\end{equation}
Expansion of~\eqref{eq:Fexpanded} gives the polynomial
in~\cite[Eq.~(1.1)]{SharipovThird}. In particular,
\begin{equation}\label{eq:constant}
P(-t)=P(t),\qquad P(0)=F(0)=(ABU)^2=(abu)^4>0.
\end{equation}
Homogeneity gives
\begin{equation}\label{eq:homogeneity}
P_{\lambda a,\lambda b,\lambda u}(\lambda t)
=\lambda^{12}P_{a,b,u}(t),\qquad \lambda\ne0.
\end{equation}

\begin{definition}[Reflection and complete degree pattern]
For a monic polynomial $h$ of degree $d$, its monic reflection is
$h^\dagger(t)=(-1)^dh(-t)$.
A complete degree pattern is the partition of $12$ formed by the
degrees of all monic irreducible factors of $P$ over $\mathbb{Q}$,
counted with multiplicities. The one-part pattern $12$ means
irreducibility and is excluded from the list of nontrivial patterns.
The notation $d^{[q]}$ means $q$ occurrences of the degree $d$.
\end{definition}

We use the standard terminology of integer partitions;
see~\cite[Chapter~1]{AndrewsPartitions}.

\begin{theorem}[Uniform normal form]\label{thm:normal}
For positive pairwise distinct rational $a,b,u$, the complete
factorization of $P=P_{a,b,u}$ over $\mathbb{Q}$ can be written as
\begin{equation}\label{eq:normal}
P(t)=\prod_{i=1}^{k}(t^2+s_i^2)\,
      \prod_{j=1}^{r}g_j(t^2)\,H(t)H^\dagger(t),
\end{equation}
where multiplicities are retained and the following statements hold:
\begin{enumerate}
\item $0\le k\le2$ and $s_i\in\mathbb{Q}^{\times}$.
\item Each $g_j$ is monic over $\mathbb{Q}$, has degree $d_j\ge2$,
and $g_j(t^2)$ is irreducible. Moreover,
\[
g_j(0)=\varepsilon_j\rho_j^2,
\qquad \varepsilon_j\in\{1,-1\},\quad
\rho_j\in\mathbb{Q}^{\times}.
\]
\item $H$ is a product containing one member of each non-even
reflected pair, with its multiplicity. Writing $m=\deg H$, one has
\begin{equation}\label{eq:degree-sign}
k+\sum_{j=1}^{r}d_j+m=6,\qquad
       \prod_{j=1}^{r}\varepsilon_j=(-1)^m.
\end{equation}
\item There is an obligatory even factor:
\begin{equation}\label{eq:marked}
\text{some }d_j\ge3\quad\text{or}\quad
\text{some }d_j=2\text{ with }\varepsilon_j=1.
\end{equation}
\end{enumerate}
In particular $r\ge1$ and $m\le4$.
\end{theorem}

\begin{corollary}[Complete numerical inventory]\label{cor:54}
The necessary conditions in Theorem~\ref{thm:normal} leave exactly
$22$ nontrivial degree patterns. Consequently $54$ of the $76$
nontrivial partitions of $12$ cannot be complete irreducible degree
patterns of $P_{a,b,u}$. Table~\ref{tab:excluded} records all $54$
exclusions, and Table~\ref{tab:remaining} lists the $22$ patterns not
excluded by these conditions.
\end{corollary}

The count describes compatibility with proved necessary conditions.
It does not assert that any remaining reducible pattern occurs for
a permitted parameter triple. Arbitrary groupings of irreducible
factors into larger reducible factors are not complete degree patterns.

\section{Root valuations, square classes, and reflection}\label{sec:global}

For $V\ne0$ write $R_V(x)=(x-V)/(x+V)$. Direct evaluation gives
\begin{align}
F(-A)&=4A^2(A-B)^2(A-U)^2,\nonumber\\
F(-B)&=4B^2(A-B)^2(B-U)^2,\nonumber\\
F(-U)&=-4U^2(A-U)^2(B-U)^2.\label{eq:special}
\end{align}
Under pairwise distinctness all these values are nonzero. Hence
every root of $F$ satisfies
\begin{equation}\label{eq:ratios}
R_A(x)^2+R_B(x)^2=R_U(x)^2.
\end{equation}

\begin{lemma}[Valuations of roots]\label{lem:valuations}
Let $a,b,u\in\mathbb{Q}^{\times}$ have pairwise distinct squares.
For every prime $p$ and every root $\theta$ of $F$ in
$\overline{\mathbb{Q}}_p$, with $v_p(p)=1$, one has
\begin{equation}\label{eq:rootvaluations}
v_p(\theta)\in\{2v_p(a),2v_p(b),2v_p(u)\}.
\end{equation}
Every monic rational divisor $g$ of $F$ therefore satisfies
\begin{equation}\label{eq:signedsquare}
g(0)\in\mathbb{Q}^{\times2}\cup-\mathbb{Q}^{\times2}.
\end{equation}
\end{lemma}
\begin{proof}
If $v_p(\theta)$ differs from all three parameter-square valuations,
each $R_V(\theta)$ is an integral unit reducing to $1$ or $-1$.
Reduction of~\eqref{eq:ratios} then gives $1+1=1$ in the residue
field, a contradiction in every characteristic, including two.
This proves~\eqref{eq:rootvaluations}.
The valuation of $g(0)$ is the sum of the valuations of its roots,
with multiplicities, so it is even at every finite prime.
A nonzero rational number having this property is a square up to sign.
\end{proof}

\begin{lemma}[Even quadratics]\label{lem:quadratics}
Every even monic irreducible quadratic factor of $P$ over
$\mathbb{Q}$ has the form $t^2+s^2$, with $s\in\mathbb{Q}^{\times}$.
\end{lemma}
\begin{proof}
If $t^2+c$ divides $P=F(t^2)$, then $-c$ is a rational root of $F$.
Lemma~\ref{lem:valuations} gives $-c=\pm s^2$.
The positive sign would split $t^2+c=t^2-s^2$ over $\mathbb{Q}$.
Thus $c=s^2$.
\end{proof}

\begin{lemma}[Reflection and descent of evenness]\label{lem:reflection}
Every non-even monic irreducible factor $h$ of $P$ occurs with its
distinct reflection $h^\dagger$, with the same multiplicity.
If a rational irreducible factor of $P$ contains a nonconstant even
factor over $\mathbb{Q}_2$, that rational irreducible factor is even.
Moreover, no non-even rational irreducible factor has a purely
imaginary root.
\end{lemma}
\begin{proof}
Reflection preserves irreducibility and multiplicity since $P$ is
even. An irreducible factor fixed by reflection cannot have odd
degree: it would be an odd polynomial divisible by $t$, contrary
to $P(0)\ne0$. Thus it is even.
For the second assertion, let $\kappa$ be an even local factor of $h$.
It also divides $h^\dagger$. If $h\ne h^\dagger$, their rational
Bezout identity persists over $\mathbb{Q}_2$ and rules out a common
nonconstant factor. Hence $h=h^\dagger$ is even.
Finally, if $h$ had a purely imaginary root $\beta$, its rational
coefficients would make $\overline\beta=-\beta$ a root as well.
Then $h$ and $h^\dagger$ would share a root, forcing
$h=h^\dagger$ by irreducibility.
\end{proof}

\section{Exactly two negative roots}\label{sec:real}

\begin{lemma}[Real root count]\label{lem:negative}
For positive pairwise distinct real $A,B,U$, the sextic $F_{A,B,U}$
has exactly two negative real roots, and both are simple.
Consequently, for the rational-square parameters in~\eqref{eq:ABU},
the polynomial $P$ has at most two even irreducible quadratic
factors, counted with multiplicities.
\end{lemma}
\begin{proof}
By symmetry take $A<B$. For $z>0$ outside $\{A,B,U\}$ put
\[
q_V(z)=\frac{R_V(z)}{R_U(z)},\qquad
L(z)=q_A(z)^{-2}+q_B(z)^{-2}.
\]
Then
\begin{align}
q_V'(z)&=\frac{2(V-U)(z^2+UV)}{(z+V)^2(z-U)^2},\label{eq:qderivative}\\
F(-z)&=(z-A)^2(z-B)^2(z+U)^2(L(z)-1).\label{eq:realidentity}
\end{align}
The function $L$ extends across $U$ with value zero. The points
$A,B,U$ themselves do not give roots of $F(-z)$,
by~\eqref{eq:special}.

If $U<A<B$, both summands of $L$ strictly decrease on $(0,U)$,
with $L$ going from $2$ to $0$, and strictly increase on $(U,A)$,
with $L$ going from $0$ to $+\infty$. For $z>A$ one has
$0<q_A(z)<1$, so $L(z)>1$.
If $A<U<B$, both summands decrease on $(A,U)$, where $L$ goes
from $+\infty$ to $0$, and increase on $(U,B)$, where it goes
from $0$ to $+\infty$. For $0<z<A$ one has $0<q_A(z)<1$,
and for $z>B$ one has $0<q_B(z)<1$, excluding further roots.
If $A<B<U$, one has $0<q_B(z)<1$ for $0<z<B$.
Both summands decrease on $(B,U)$ from a total limit $+\infty$
to $0$, and increase on $(U,\infty)$ from $0$ to $2$.

These signs follow directly from~\eqref{eq:qderivative} and
$(q_V^{-2})'=-2q_V'/q_V^3$. There are precisely two crossings of
$L=1$ in each ordering. At a crossing the two summand derivatives
have the same nonzero sign. Thus the crossings, and hence the
negative roots of $F$, are simple.

By Lemma~\ref{lem:quadratics}, an even irreducible quadratic of $P$
corresponds to a negative rational root $-s^2$ of $F$. There are at
most two such roots, and they are simple, proving the final assertion.
\end{proof}

\section{The local alternative and its tools}\label{sec:local}

In the local argument $a,b,u\in\mathbb{Q}_2^\times$, and their
squares are pairwise distinct. Write $v=v_2$, with $v(2)=1$,
extended to the algebraic closure, and put $K=\mathbb{Q}_2(i)$.
We use the elementary Newton polygon and Hensel principles for local
fields~\cite{Serre,Robert}. Weierstrass preparation over $\mathbb{Z}_2$
is used to isolate roots with positive valuation; the formal version
needed here is given in~\cite[Corollary~1.2]{Berger}.

A Newton edge of length $n$ and slope $-h$ supplies $n$ roots of
valuation $h$, counted with multiplicities. The reduced denominator
of $h$ divides the ramification index of the field generated by
any such root. Root clusters separated by valuation give factors
over the ground field. For an integral power series with first
unit coefficient in degree $n$, Weierstrass preparation gives an
integral power-series unit times a monic distinguished polynomial
of degree $n$. Its roots are the positive-valuation roots under
consideration. In each use below the series is rational with
denominator an analytic unit, so no root is lost on its domain.

\begin{lemma}[Square class and degree doubling]\label{lem:doubling}
Let $M/\mathbb{Q}_2$ be finite. If $y\in M$ and $v(y)>2$, then
$1-y$ is a square in $M$. In particular, if
\[
\theta=-u^2+q,\qquad M=\mathbb{Q}_2(q),\qquad v(q/u^2)>2,
\]
then $\theta=-u^2s^2$ for some $s\in M^\times$.
If $i\notin M$ and $d=[M:\mathbb{Q}_2]$, the minimal polynomial
of a square root of $\theta$ is even of degree $2d$ and has square
constant term in $\mathbb{Q}_2$.
\end{lemma}
\begin{proof}
Hensel's lemma applied at $1$ to $Z^2-(1-y)$ gives the first
assertion because the error has valuation greater than twice the
derivative valuation; see~\cite[Chapter~II, \S1.5]{Robert}. Let $f$ be the minimal polynomial of $\theta$.
The fields $\mathbb{Q}_2(\theta)$ and $\mathbb{Q}_2(q)$ agree.
If $\beta^2=\theta=-u^2s^2$, then
$\mathbb{Q}_2(\beta)=M(i)$ has degree $2d$.
Hence $f(t^2)$ is its minimal polynomial. Its constant term is
\begin{equation}\label{eq:normconstant}
f(0)=(-1)^dN_{M/\mathbb{Q}_2}(\theta)
=u^{2d}N_{M/\mathbb{Q}_2}(s)^2.
\end{equation}
\end{proof}

\begin{lemma}[Two facts about $\mathbb{Q}_2(i)$]\label{lem:K}
The field $K=\mathbb{Q}_2(i)$ has value group
$\tfrac12\mathbb{Z}$, residue field $\mathbb{F}_2$, and ring of
integers $\mathbb{Z}_2[i]$. Every unit $z$ in this ring satisfies
\begin{equation}\label{eq:fourthpower}
z^4\equiv1\pmod{8\mathcal O_K}.
\end{equation}
Furthermore $-1$ is not a norm from $K$ to $\mathbb{Q}_2$.
\end{lemma}
\begin{proof}
The polynomial of $1+i$ is $X^2-2X+2$, which is Eisenstein;
see also~\cite[Chapter~II, \S4.5]{Robert} for the Gaussian
$2$-adic extension.
For a unit $z=r+is$, $r,s\in\mathbb{Z}_2$, exactly one of $r,s$
is odd. The real part of $z^2$ is odd and its imaginary part is
divisible by $4$, giving~\eqref{eq:fourthpower} after squaring.
If $r^2+s^2=-1$ with a negative valuation among $r,s$, unequal
valuations give a uniquely minimal square, and equal valuations
$m<0$ give valuation $2m+1<0$ for their sum. Both are impossible.
For integral $r,s$ reduction modulo four is impossible.
Thus $-1$ is not a norm. Multiplication by a nonzero rational square
does not change membership in the norm group.
\end{proof}

\begin{theorem}[Uniform local alternative]\label{thm:local}
For nonzero $a,b,u\in\mathbb{Q}_2$ with pairwise distinct squares,
$P_{a,b,u}$ has either an even irreducible factor of degree six or
eight, or an even monic irreducible quartic with square constant term
in $\mathbb{Q}_2$. In the quartic case supplied by the construction,
the root field is biquadratic over $\mathbb{Q}_2$.
\end{theorem}

The proof occupies Sections~\ref{sec:unique}--\ref{sec:units}.
By~\eqref{eq:homogeneity}, dividing all parameters and $t$ by a
common power of two preserves factor degrees and evenness.
A quartic constant changes by a fourth power, preserving its
square class.

\section{A unique distinguished valuation}\label{sec:unique}

\begin{lemma}\label{lem:unique}
If $v(u)$ differs from both $v(a)$ and $v(b)$, then $P$ has an even
irreducible quartic over $\mathbb{Q}_2$ with root field
$\mathbb{Q}_2(i,\sqrt2)$ and square constant term.
\end{lemma}
\begin{proof}
Put $t=u(1+z)$. For $w=a,b$, choose
\[
c_w=\begin{cases}
w/u,&v(w)>v(u),\\
u/w,&v(w)<v(u).
\end{cases}
\qquad e_w=c_w^2/4\in\mathbb{Z}_2.
\]
Choose the integral numerator and denominator
\[
(N_w,D_w)=\begin{cases}
((1+z)^2-c_w^2,(1+z)^2+c_w^2),&v(w)>v(u),\\
(c_w^2(1+z)^2-1,c_w^2(1+z)^2+1),&v(w)<v(u).
\end{cases}
\]
Then $R_{w^2}(t^2)=N_w/D_w$, $D_w(0)$ is a unit, and
$N_w^2=D_w^2-16e_w(1+z)^2$.
With
\[
D=D_a^2D_b^2,\qquad
T=(1+z)^2(e_aD_b^2+e_bD_a^2),
\]
the root equation becomes
\begin{equation}\label{eq:uniqueG}
G(z)=D(z)Q(z)+16T(z)(z^2+2z+2)^2=0,
\end{equation}
where
\[
Q(z)=z^2(z+2)^2-2(z^2+2z+2)^2
=-(z^4+4z^3+12z^2+16z+8).
\]
Writing $G=\sum g_jz^j$, one obtains
\[
v(g_0)=3,\quad v(g_1)\ge3,\quad v(g_2)\ge2,
\quad v(g_3)\ge2,\quad v(g_4)=0,
\]
and all higher coefficients are integral. For example, the constant
correction is divisible by $64$, and only $-D(0)$ contributes an odd
degree-four coefficient. Thus the Newton polygon has an edge from
$(0,3)$ to $(4,0)$, supplying exactly four roots with $v(z)=3/4$.
All cleared denominators are nonzero there: each $D_w$ is a unit,
and $v(z^2+2z+2)=1$.

Each root has degree divisible by four. Since there are exactly four
roots of this valuation, they constitute one irreducible quartic
block. Reflection sends $z$ to $-2-z$ and preserves this cluster,
so its polynomial in $t$ is even.

For a root $\beta=t$ put $\theta=\beta^2$ and
$M=\mathbb{Q}_2(\theta)$. The irreducible even quartic gives
$[M:\mathbb{Q}_2]=2$. Set $R=(\theta-U)/(\theta+U)$ and
$w=T(z)/D(z)$. Although written in $z$, this last expression belongs
to $M$, since
$16w=2-R_A(\theta)^2-R_B(\theta)^2$; it is integral.
At a root of~\eqref{eq:uniqueG},
\[
R^2=2(1-8w),\qquad R^2-1=1-16w.
\]
Both $1-8w$ and $1-16w$ are squares in $M$ by Hensel's lemma.
Hence $M=\mathbb{Q}_2(\sqrt2)$, and
\[
\frac{(1+R)^2}{R^2-1}=-\frac{\theta}{U}
\]
gives $\theta=-u^2s^2$ in $M$.
The field $M$ does not contain $i$, so the root field of the quartic
is $M(i)=\mathbb{Q}_2(i,\sqrt2)$. Its constant term is a square
by the norm calculation~\eqref{eq:normconstant}.
\end{proof}

\section{Exactly two parameters share the distinguished valuation}\label{sec:shared}

\begin{lemma}\label{lem:shared}
Suppose $v(a)=v(u)\ne v(b)$ and $a^2\ne u^2$.
After scaling $a,u$ to units put $r=v(a^2-u^2)\ge3$.
The polynomial $P$ has an even irreducible local factor of degree
eight if $r=3$, degree six if $r>3$ and $3\nmid r$, and degree
four with square constant term if $r>3$ and $3\mid r$.
In the last case its root field is biquadratic.
\end{lemma}
\begin{proof}
Take $a,u$ to be units, write $C=A-U$, and set $q=x+U$.
Since $v(B)$ is a nonzero even integer,
\[
R_B(-U+q)^2=1+16f(q),\qquad f(q)\in\mathbb{Z}_2[[q]].
\]
For $v(B)>0$ this follows from
$R_B(x)^2=1-4Bx/(x+B)^2$; for $v(B)<0$ use the same expression
after dividing numerator and denominator by $B$.
The denominators are analytic units on $v(q)>0$.
The root equation is equivalent there to
\begin{equation}\label{eq:sharedg}
g(q)+16f(q)q^2(q+C)^2=0,
\end{equation}
where
\begin{equation}\label{eq:quarticg}
g(q)=q^4-2Cq^3+C(C+12U)q^2
 +4UC(C-2U)q-4U^2C^2.
\end{equation}
The valuations of its coefficients in increasing degree are
\[
2r+2,\quad r+3,\quad r+2,\quad r+1,\quad0.
\]
The perturbation in~\eqref{eq:sharedg} is above the relevant lower
edges. The full integral power series reduces to $q^4$ modulo two;
Weierstrass preparation therefore supplies a degree-four factor
capturing exactly its positive-valuation roots.
The selected roots below satisfy $q\ne0,-C$,
and $x+B\ne0$, so no cleared pole produces a root.

If $r=3$, the Newton polygon is the segment $(0,8)$--$(4,0)$.
Set $q=4z$ and divide by $2^8$. The normalized equation is integral
and its reduction is
\[
z^4+z+1\in\mathbb{F}_2[z].
\]
It has no linear root and is not divisible by $z^2+z+1$, so it is
irreducible. Hensel lifting gives a root in the unramified quartic
extension $L$, and the residue of $z$ generates its residue field;
therefore $\mathbb{Q}_2(q)=L$.
For $\theta=q-U$, one has $\theta/U\equiv-1\pmod4$ in $L$.
If $\theta/U$ were a square, its square root would have residue $1$,
and hence be $1+2h$ because $L$ is unramified. Its square would be
$1$ modulo four, a contradiction. Thus adjoining $\sqrt\theta$
doubles the degree to eight, and its minimal polynomial is the even
polynomial obtained by composing that of $\theta$ with $t^2$.

If $r>3$, the lower vertices are
\[
(0,2r+2),\quad(1,r+3),\quad(4,0).
\]
When $3\nmid r$, the last segment has length three and root valuation
$(r+3)/3$ of denominator three. It gives an irreducible cubic factor
and a field $L$ of degree three. Its selected roots have valuation
greater than two, and $i\notin L$. Lemma~\ref{lem:doubling} gives
an even degree-six factor of $P$.

Finally let $r=3m>3$, so $m\ge2$. Substitute $q=2^{m+1}z$ and
divide~\eqref{eq:sharedg} by $2^{4m+4}$. All coefficients are
integral. The coefficients of degrees zero through four have
valuations $2m-2,0,m,2m,0$, respectively; all coefficients of the
normalized perturbation are divisible by two. The reduction is
\[
z^4+z=z(z+1)(z^2+z+1).
\]
The simple irreducible quadratic factor lifts in the unramified
quadratic extension $L$; its nontrivial residue guarantees that
$\mathbb{Q}_2(q)=L$. The substitutions in the power series converge
on integral $z$, so Hensel lifting is applicable.
Now $v(q)=m+1\ge3$ and $i\notin L$. Lemma~\ref{lem:doubling}
gives an even quartic with square constant term and root field $L(i)$,
which is biquadratic.
\end{proof}

\section{Three equal parameter valuations}\label{sec:units}

Scale the parameters to units. Then $A,B,U\equiv1\pmod8$.
Put
\[
C=A-U,\qquad D=B-U,\qquad c=v(C),\quad d=v(D),
\]
where $c,d\ge3$. By symmetry assume $c\le d$.
For $q=x+U$, write
\[
\frac{F(-U+q)}{4U^2}=\sum_{j=0}^6a_jq^j.
\]
Exact expansion gives
\begin{align}
a_6&=1/(4U^2),&
a_5&=-(C+D+2U)/(2U^2),\nonumber\\
a_4&=1+(C+D)/U+(C^2-12CD+D^2)/(4U^2),\nonumber\\
a_3&=-CD(C+D-16U)/(2U^2),\nonumber\\
a_2&=CD(CD+12U(C+D)-16U^2)/(4U^2),\nonumber\\
a_1&=CD(CD-2U(C+D))/U,&
a_0&=-C^2D^2.\label{eq:unitcoeffs}
\end{align}
In particular,
\begin{gather}
v(a_6)=-2,\quad v(a_5)=0,\quad a_4\equiv1\pmod8,\nonumber\\
v(a_3)\ge c+d+2,\quad v(a_2)=c+d+2,\label{eq:unitvals}\\
v(a_1)=2c+d+1\quad(c<d),\qquad
v(a_1)\ge3c+2\quad(c=d).\nonumber
\end{gather}
The rightmost Newton edge is $(4,0)$--$(6,-2)$; all earlier
points lie on or above the line of slope $-3$ through $(4,0)$.
Thus there are two roots with $v(q)=1$ and four with $v(q)\ge3$.
Let their monic cluster polynomials be $V,W$, so that
\begin{equation}\label{eq:VW}
F(-U+q)=V(q)W(q),\qquad \deg V=2,\quad\deg W=4.
\end{equation}

\begin{lemma}[The four-root cluster does not split in $K$]\label{lem:nonsplit}
The four roots of $W$ do not all belong to $K=\mathbb{Q}_2(i)$.
\end{lemma}
\begin{proof}
Let $g=d-c$. From~\eqref{eq:unitcoeffs}--\eqref{eq:unitvals},
the valuations of the four roots are
\begin{equation}\label{eq:gapvals}
\begin{array}{c|c}
g&\text{root valuations}\\\hline
0&c,c,c,c\\
1&c+\tfrac12,c+\tfrac12,c+\tfrac12,c+\tfrac12\\
2&c+1,c+1,c+1,c+1\\
g>2&d-1,\ c+(g+1)/3,\ c+(g+1)/3,\ c+(g+1)/3.
\end{array}
\end{equation}
We consider these four possibilities.

If $g=1$, write $C=2^cC_0$, $D=2^{c+1}D_0$, $e=C_0D_0$.
A hypothetical root in $K$ has the form $q=2^c(1+i)z$ with $z$ a
unit. By~\eqref{eq:fourthpower},
\[
v(q^4-C^2D^2)=v(-2^{4c+2}(z^4+e^2))=4c+3.
\]
But $v(a_1q)=4c+5/2$. All remaining terms, after grouping
$q^4+a_0$, have larger valuations, including $(a_4-1)q^4$.
The minimum is unique, so the equation cannot hold.

If $g=2$, substitute $q=2^{c+1}z$ and divide by $2^{4c+4}$.
The reduction is $z^4+z+1$. A unit of $K$ reduces to $1$, which
is not a root of that polynomial. None of these four roots is in $K$.

If $g>2$, put $h=c+(g+1)/3$. If $h$ is not integral, its
denominator is three, so $h\notin\tfrac12\mathbb{Z}$; the three
roots of valuation $h$ cannot be in $K$. If $h$ is integral,
substitute $q=2^hz$ and divide by $2^{4h}$. The resulting integral
polynomial reduces to $z^4+z$. The three roots of valuation $h$
give units $z$, but the only nonzero residue in $\mathbb{F}_2$ is
$1$, a simple root. Hensel uniqueness permits only one such root
in $K$, so these three roots cannot all belong to $K$.

It remains to treat $g=0$. Put
\[
C=2^cC_0,\quad D=2^cD_0,\quad
e=C_0D_0,\quad j=(C_0+D_0)/2,\quad h=2^c,
\]
and
\[
T(z)=\frac{F(-U+hz)}{4U^2h^4}.
\]
There is a useful direct identity
\begin{align}
T(z)={}&\left[z(z+D_0)\left(1+\frac{h(C_0-z)}{2U}\right)\right]^2
\nonumber\\
&+\left[z(z+C_0)\left(1+\frac{h(D_0-z)}{2U}\right)\right]^2
\nonumber\\
&-\left[(z+C_0)(z+D_0)\left(1-\frac{hz}{2U}\right)\right]^2.
\label{eq:Tidentity}
\end{align}
Since $v(h)\ge3$, each additional factor in parentheses has square
congruent to $1$ modulo eight. Thus
\begin{equation}\label{eq:Tcongruence}
T(z)\equiv z^4-e^2-4ez^2-4ejz
\pmod{8\mathcal O_K[z]}.
\end{equation}
For a unit root $z\in K$, equation~\eqref{eq:fourthpower} and
$e^2\equiv1\pmod8$ imply $z+j\in2\mathcal O_K$.
If $j$ is even, this contradicts that $z$ is a unit.
If $j$ is odd, it implies $v(z-1)\ge1$.
On the other hand~\eqref{eq:Tcongruence} gives
\[
T(1+w)\equiv w^4+2w^2\pmod{4\mathcal O_K[w]}.
\]
The coefficient valuations, in increasing degree, satisfy
\[
(\ge2,\ \ge2,\ 1,\ \ge2,\ 0,\ \ge3,\ \ge4).
\]
Thus every coefficient point $(\ell,v_2([w^\ell]T(1+w)))$ is on
or above the supporting line $y=2-\ell/2$, which passes through
$(2,1)$ and $(4,0)$.
The Newton polygon therefore has an edge of
slope $-1/2$ through $(2,1)$ and $(4,0)$, possibly extended left
to $(0,2)$. There are at least two roots with $v(w)=1/2$.
They give unit roots $z=1+w$, hence roots of the four-root cluster,
but cannot be in $K$ because roots there required $v(z-1)\ge1$.
This completes all cases.
\end{proof}

\begin{lemma}[The norm of the four-root cluster]\label{lem:Wnorm}
For the factorization~\eqref{eq:VW},
\begin{equation}\label{eq:Wnorm}
V(0)\in\mathbb{Q}_2^{\times2},\qquad
W(0)\in-\mathbb{Q}_2^{\times2}.
\end{equation}
There is a root $q$ of $W$ whose field $M=\mathbb{Q}_2(q)$ has
degree $d_0\in\{2,3,4\}$ and does not contain $i$.
\end{lemma}
\begin{proof}
Define
\[
\mathcal H(Y)=\frac{F(U+4Y)}{256U^4}.
\]
The original three-square identity gives
\begin{align}
\mathcal H(Y)={}&\left[(Y-C/4)
 \left(1+\frac{D+4Y}{2U}\right)
 \left(1+\frac{2Y}{U}\right)\right]^2\nonumber\\
&+\left[(Y-D/4)
 \left(1+\frac{C+4Y}{2U}\right)
 \left(1+\frac{2Y}{U}\right)\right]^2\nonumber\\
&-\left[Y
 \left(1+\frac{C+4Y}{2U}\right)
 \left(1+\frac{D+4Y}{2U}\right)\right]^2.
\label{eq:Hidentity}
\end{align}
Every bracket is an integral polynomial congruent to $Y$ modulo two.
Its square is therefore congruent to $Y^2$ modulo four, and
\begin{equation}\label{eq:Hcongruence}
\mathcal H(Y)\equiv Y^2\pmod{4\mathbb{Z}_2[Y]}.
\end{equation}
Weierstrass preparation gives
\[
\mathcal H(Y)=J(Y)(Y^2+\beta Y+e),
\]
where $J$ is an integral power-series unit and initially
$\beta,e\in2\mathbb{Z}_2$. The constant coefficient
in~\eqref{eq:Hcongruence} gives $e\in4\mathbb{Z}_2$; comparison of
linear coefficients then gives $\beta\in4\mathbb{Z}_2$.
The two distinguished roots $Y_1,Y_2$ have positive valuation.
They give the two roots $q_j=2U+4Y_j$ of valuation one, and hence
\[
V(0)=q_1q_2=4U^2-8U\beta+16e
=4U^2(1-2\beta/U+4e/U^2).
\]
The parenthesis is $1$ modulo eight, so it is a square.
Using~\eqref{eq:special} and~\eqref{eq:VW},
\[
W(0)=\frac{F(-U)}{V(0)}
=\frac{-4U^2C^2D^2}{V(0)}\in-\mathbb{Q}_2^{\times2}.
\]

Choose a root $q\notin K$ using Lemma~\ref{lem:nonsplit}, and put
$M=\mathbb{Q}_2(q)$, $d_0=[M:\mathbb{Q}_2]\le4$.
Then $d_0\ge2$. For $d_0=2$, the fact that $q\notin K$ excludes
$M=K$; for $d_0=3$, degree parity excludes $K\subset M$.
If $d_0=4$ and $K\subset M$, then $W$ is the minimal polynomial
of $q$ and
\[
W(0)=N_{M/\mathbb{Q}_2}(q)
=N_{K/\mathbb{Q}_2}(N_{M/K}(q)).
\]
This is impossible by~\eqref{eq:Wnorm} and Lemma~\ref{lem:K}.
Thus $i\notin M$ in every case.
\end{proof}

\begin{proof}[Proof of Theorem~\ref{thm:local}]
If $v(u)$ differs from both other parameter valuations, apply
Lemma~\ref{lem:unique}. If it agrees with exactly one, interchange
$a,b$ if necessary and use Lemma~\ref{lem:shared}.
If all three agree, scale to units and choose $q,M,d_0$ as in
Lemma~\ref{lem:Wnorm}. Since $v(q)\ge3$ and $U$ is a unit,
Lemma~\ref{lem:doubling} applied to $\theta=-U+q$ gives an even
irreducible factor of degree $2d_0$ with square constant term.
For $d_0=3,4$ this gives degree six or eight. For $d_0=2$, its
root field $M(i)$ is the compositum of two distinct quadratic
extensions, hence biquadratic. The three valuation cases exhaust
all parameter triples.
\end{proof}

\section{Globalization and proof of the normal form}\label{sec:normalproof}

\begin{proposition}[The obligatory global factor]\label{prop:obligatory}
For positive pairwise distinct rational $a,b,u$, the polynomial $P$
has an even irreducible factor of degree at least six or an even
irreducible quartic $t^4+pt^2+\rho^2$, with
$p\in\mathbb{Q}$ and $\rho\in\mathbb{Q}^{\times}$.
\end{proposition}
\begin{proof}
Choose a local factor $\kappa$ from Theorem~\ref{thm:local}, and
a monic rational irreducible factor $h$ of $P$ divisible by
$\kappa$ over $\mathbb{Q}_2$.
Lemma~\ref{lem:reflection} forces $h$ to be even.
If $\deg\kappa\ge6$ or $\deg h>4$, then $\deg h\ge6$.
Otherwise $\deg h=\deg\kappa=4$, so $h=\kappa$ is an even
quartic with square constant term in $\mathbb{Q}_2$.
Write $h(t)=g(t^2)$, $g=x^2+px+c$.
Polynomial division shows that $g\mid F$: a remainder of degree
less than two would, after composition, be divisible by $h$ and
have degree less than four. Lemma~\ref{lem:valuations} gives
$c=\pm\rho^2$ over $\mathbb{Q}$.
The negative sign is incompatible with $c\in\mathbb{Q}_2^{\times2}$
because $-1$ is not a square in $\mathbb{Q}_2$.
Thus $c=\rho^2>0$.
\end{proof}

\begin{proof}[Proof of Theorem~\ref{thm:normal}]
Separate the monic rational irreducible factors fixed by reflection
from the distinct reflected pairs. Choose one member from each
pair, with multiplicities, and let $H$ be their product.
The other members have product $H^\dagger$.
The fixed irreducibles are even by Lemma~\ref{lem:reflection}.
Their quadratic members have the form $t^2+s_i^2$ by
Lemma~\ref{lem:quadratics}, and their number $k$ is at most two
by Lemma~\ref{lem:negative}.

Write every remaining even irreducible as $g_j(t^2)$, with $g_j$
monic of degree $d_j\ge2$. The same polynomial-division argument
as above shows $g_j\mid F$; it is irreducible since a factorization
of $g_j$ would factor $g_j(t^2)$. Lemma~\ref{lem:valuations} gives
$g_j(0)=\varepsilon_j\rho_j^2$.
This proves~\eqref{eq:normal}. Degree comparison gives the first
equation in~\eqref{eq:degree-sign}.

Since
\[
H(0)H^\dagger(0)=(-1)^mH(0)^2
\]
and $P(0)>0$, comparison of constant terms gives the second equation.
Finally Proposition~\ref{prop:obligatory} gives~\eqref{eq:marked}.
It forces $\sum d_j\ge2$, hence $r\ge1$ and $m\le4$.
\end{proof}

\begin{corollary}[Degree and parity consequences]\label{cor:parity}
Every non-even irreducible rational factor of $P$ has degree at
most four. In particular no irreducible odd degree above three
occurs. If all factors have degree at most four and exactly one
quartic occurs, then the number of quadratic factors is even.
\end{corollary}
\begin{proof}
The first statements follow from $m\le4$, since each non-even
factor appears in $H$ or $H^\dagger$, and an even irreducible
cannot have odd degree.
Under the last hypothesis, the unique quartic is even by reflection.
It is the only factor $g_j(t^2)$ in~\eqref{eq:normal}; thus $r=1$,
$d_1=2$, and $\varepsilon_1=1$ by~\eqref{eq:marked}.
The sign equation makes $m$ even, and $k+2+m=6$ makes $k$ even.
Every other quadratic occurs in a reflected pair. Hence the total
number of quadratics is even.
\end{proof}

\begin{corollary}[A common sign obstruction]\label{cor:sign}
The three complete patterns
\[
4+3+3+2,\qquad4+2+2+2+1+1,\qquad4+2+1^{[6]}
\]
are impossible for every positive pairwise distinct rational triple.
\end{corollary}
\begin{proof}
Each has a unique quartic, no larger factor, and an odd number of
quadratics, contrary to Corollary~\ref{cor:parity}.
\end{proof}

\begin{remark}[The role of all finite primes]
A rational number that is a square in $\mathbb{Q}_2$ need not be
positive; for example $-7$ is a square there.
The positivity in Proposition~\ref{prop:obligatory} uses both the
global signed-square statement~\eqref{eq:signedsquare}, proved at
every finite prime, and the local square class at two.
Similarly, a local degree partition is not automatically a rational
degree partition. Reflection and the persistent Bezout identity
are what allow the local even factor to control a global irreducible.
\end{remark}

\section{Excluded complete degree patterns}\label{sec:excluded}

We now translate the proved factor restrictions into a complete
inventory of exclusions. Throughout this inventory, a pattern records
all irreducible factors over $\mathbb{Q}$, with multiplicities, as in
the definition in Section~\ref{sec:poly}. Each exclusion holds for every
positive pairwise distinct rational triple. The three steps below are
applied successively, so that a pattern already excluded at one step
is not counted again at a later step.

\medskip
\noindent\textbf{(I) Reflection.}
By Lemma~\ref{lem:reflection}, every irreducible factor of odd degree
occurs together with a distinct reflected factor of the same degree
and multiplicity. Consequently each odd degree must occur an even
number of times in a complete pattern. Among the $76$ nontrivial
partitions of twelve, exactly $37$ fail this condition. They are
already excluded by evenness, before any additional arithmetic input
is used.

\medskip
\noindent\textbf{(II) The obligatory even factor.}
Of the $39$ patterns that pass the reflection test, a further $14$
contain no even degree of size at least four.
Proposition~\ref{prop:obligatory} excludes all of them: every complete
factorization must contain an even irreducible factor of degree at
least six or an even irreducible quartic with positive square constant
term. In particular, the degree pattern must contain an even degree
at least four. Passing this degree test does not yet establish the
required evenness or constant term of a particular factor; it only
identifies the patterns that this step has not excluded.

\medskip
\noindent\textbf{(III) The sign condition for a unique quartic.}
There are now $25$ patterns left. Corollary~\ref{cor:parity} requires
an even number of quadratic factors whenever there is exactly one
quartic and no larger factor. Corollary~\ref{cor:sign} explicitly
excludes the three patterns violating this requirement.
This last step uses the sign equation and the positive-square constant
of the obligatory quartic, rather than reflection alone.

Thus the successive exclusions number $37+14+3=54$, leaving
$76-54=22$ nontrivial degree patterns to examine for compatibility
with the normal form. These are successive counts of abstract degree
patterns, not independent Diophantine obstructions.
Table~\ref{tab:excluded} records the full proved exclusion status:
every displayed pattern is impossible, and the heading of each group
identifies the result that excludes it. Within each group, patterns
are arranged by their largest degree.

\begingroup
\small
\setlength{\LTleft}{0pt}
\setlength{\LTright}{0pt}
\renewcommand{\arraystretch}{1.18}
\begin{longtable}{@{}p{0.90\textwidth}r@{}}
\caption{All $54$ strictly excluded complete irreducible degree patterns.
The three groups are successive and disjoint.}\label{tab:excluded}\\
\toprule
Strictly excluded complete patterns & Count\\
\midrule
\endfirsthead
\multicolumn{2}{@{}l}{\tablename~\thetable\ (continued)}\\
\toprule
Strictly excluded complete patterns & Count\\
\midrule
\endhead
\midrule
\multicolumn{2}{r@{}}{Continued on the next page}\\
\endfoot
\bottomrule
\endlastfoot
\multicolumn{2}{@{}l}{\textit{(I) Reflection: Lemma~\ref{lem:reflection}}}\\*
$11+1$,\quad $9+3$,\quad $9+2+1$\newline $9+1^{[3]}$,\quad $8+3+1$,\quad $7+5$\newline $7+4+1$,\quad $7+3+2$,\quad $7+3+1^{[2]}$\newline $7+2^{[2]}+1$,\quad $7+2+1^{[3]}$,\quad $7+1^{[5]}$\newline $6+5+1$,\quad $6+3+2+1$,\quad $6+3+1^{[3]}$\newline $5+4+3$,\quad $5+4+2+1$,\quad $5+4+1^{[3]}$\newline $5+3^{[2]}+1$,\quad $5+3+2^{[2]}$,\quad $5+3+2+1^{[2]}$\newline $5+3+1^{[4]}$,\quad $5+2^{[3]}+1$,\quad $5+2^{[2]}+1^{[3]}$\newline $5+2+1^{[5]}$,\quad $5+1^{[7]}$,\quad $4^{[2]}+3+1$\newline $4+3+2^{[2]}+1$,\quad $4+3+2+1^{[3]}$,\quad $4+3+1^{[5]}$\newline $3^{[3]}+2+1$,\quad $3^{[3]}+1^{[3]}$,\quad $3+2^{[4]}+1$\newline $3+2^{[3]}+1^{[3]}$,\quad $3+2^{[2]}+1^{[5]}$,\quad $3+2+1^{[7]}$\newline $3+1^{[9]}$ & $37$\\
\midrule
\multicolumn{2}{@{}l}{\textit{(II) Obligatory even factor: Proposition~\ref{prop:obligatory}}}\\*
$5^{[2]}+2$,\quad $5^{[2]}+1^{[2]}$,\quad $3^{[4]}$\newline $3^{[2]}+2^{[3]}$,\quad $3^{[2]}+2^{[2]}+1^{[2]}$,\quad $3^{[2]}+2+1^{[4]}$\newline $3^{[2]}+1^{[6]}$,\quad $2^{[6]}$,\quad $2^{[5]}+1^{[2]}$\newline $2^{[4]}+1^{[4]}$,\quad $2^{[3]}+1^{[6]}$,\quad $2^{[2]}+1^{[8]}$\newline $2+1^{[10]}$,\quad $1^{[12]}$ & $14$\\
\midrule
\multicolumn{2}{@{}l}{\textit{(III) Unique-quartic sign: Corollary~\ref{cor:sign}}}\\*
$4+3^{[2]}+2$,\quad $4+2^{[3]}+1^{[2]}$,\quad $4+2+1^{[6]}$ & $3$\\
\midrule
Total strictly excluded & $54$\\
\end{longtable}
\endgroup

\section{The twenty-two remaining degree patterns}\label{sec:inventory}

After all exclusions in Table~\ref{tab:excluded}, exactly $22$
nontrivial partitions remain. Table~\ref{tab:remaining} lists them
separately to distinguish the proved impossibility of the first table
from the fact that the patterns in the second are not excluded.
Each entry is compatible
with the necessary conditions of Theorem~\ref{thm:normal}; this is a
statement about abstract factor data, not the existence of a permitted
parameter triple realizing the pattern. The enumeration and this
formal compatibility are verified immediately after the table.

\begingroup
\small
\renewcommand{\arraystretch}{1.12}
\begin{longtable}{clc}
\caption{The $22$ nontrivial complete degree patterns not excluded by
Theorem~\ref{thm:normal}. Their occurrence is not asserted.}\label{tab:remaining}\\
\toprule
Largest degree & Complete patterns not excluded & Count\\
\midrule
$10$ & $10+2,\quad10+1^{[2]}$ & $2$\\[3pt]
$8$ & $8+4,\quad8+2^{[2]}$ & \\
& $8+2+1^{[2]},\quad8+1^{[4]}$ & $4$\\[3pt]
$6$ & $6^{[2]},\quad6+4+2,\quad6+4+1^{[2]}$ & \\
& $6+3^{[2]},\quad6+2^{[3]}$ & \\
& $6+2^{[2]}+1^{[2]},\quad6+2+1^{[4]},\quad6+1^{[6]}$ & $8$\\[3pt]
$4$ & $4^{[3]},\quad4^{[2]}+2^{[2]},\quad4^{[2]}+2+1^{[2]}$ & \\
& $4^{[2]}+1^{[4]},\quad4+2^{[4]}$ & \\
& $4+2^{[2]}+1^{[4]},\quad4+3^{[2]}+1^{[2]},\quad4+1^{[8]}$ & $8$\\
\bottomrule
\end{longtable}
\endgroup

\begin{proof}[Proof of Corollary~\ref{cor:54}]
The normal form provides a direct finite enumeration.
Choose $k\in\{0,1,2\}$, a nonempty list of degrees $d_j\ge2$,
and a partition $(h_1,\ldots,h_s)$ of
\[
m=6-k-\sum d_j\ge0.
\]
The corresponding degree pattern is the multiset
\begin{equation}\label{eq:projection}
\underbrace{2,\ldots,2}_{k\text{ times}},\qquad
2d_1,\ldots,2d_r,\qquad
h_1,h_1,\ldots,h_s,h_s.
\end{equation}
Keep only lists admitting signs satisfying~\eqref{eq:degree-sign}
and~\eqref{eq:marked}; remove the one-part pattern $12$.
If some $d_j\ge3$, the obligatory-factor condition imposes no sign
restriction. If all $d_j=2$ and $r\ge2$, their signs can have
either prescribed product while at least one sign is positive.
If $r=1$ and $d_1=2$, its sign must be positive and $m$ must be even.
These observations completely describe the finite enumeration.

For an explicit count, let $M$ be the largest degree in a remaining
nontrivial pattern. Corollary~\ref{cor:parity} and the obligatory
factor give $M\in\{4,6,8,10\}$.
For $M=10$, the remaining degree two is $2$ or $1+1$.
For $M=8$, reflection leaves the four possibilities
$4$, $2+2$, $2+1+1$, and $1^{[4]}$.
For $M=6$, there may be a second sextic; otherwise the remaining
degree six has paired odd degrees and gives seven possibilities.

For $M=4$, three quartics give one pattern, and two quartics leave
$2+2$, $2+1+1$, or $1^{[4]}$.
For a unique quartic, let $n_j$ be the number of degree-$j$ factors.
Then $n_1,n_3$ are even by reflection and $n_2$ is even by
Corollary~\ref{cor:parity}. Solving
$3n_3+2n_2+n_1=8$ gives four possibilities:
$3+3+1+1$, $2^{[4]}$, $2^{[2]}+1^{[4]}$, or $1^{[8]}$.
Together, these cases give the eight entries with largest degree four
listed in Table~\ref{tab:remaining}.

For completeness, all listed patterns are compatible with the
abstract normal-form constraints. For largest degree at least six,
take the displayed factors of degree at least four to be even and
pair all cubics and linears. If the number of quadratics is odd,
take one to be even and pair the others; if it is even, pair them all.
The signs can be chosen as above because a $d_j\ge3$ is present.
For largest degree four, the three patterns with multiple quartics
and small remaining degrees, together with $4^{[3]}$, admit the
same choice and appropriate signs. For a unique quartic take its
sign positive; in each listed case choose even $k$ and pair the
remaining quadratics. This satisfies both equations and the
obligatory-factor condition. These are formal compatibility choices,
not constructions of parameter triples.

There are $77$ partitions of twelve, including the one-part partition
$12$, and hence $76$ nontrivial partitions. Table~\ref{tab:remaining} has
$2+4+8+8=22$ distinct entries. Together with the $54$ disjoint exclusions in
Table~\ref{tab:excluded}, these entries account for all $76$
nontrivial partitions. Since each of the $22$ entries admits the
formal choices just described, none is excluded by the stated
normal-form conditions. This proves the exact numerical inventory
without asserting that any remaining pattern occurs for the polynomial
family. Appendix~\ref{app:certificate} independently reproduces the
enumeration from~\eqref{eq:projection} and verifies the completeness
of the underlying partition list.
\end{proof}

\section{Scope of the result}\label{sec:scope}

The normal form is a uniform necessary condition for complete
rational factorization. Its proof uses the same arithmetic mechanism
for all parameter triples: root valuations at finite primes, a real
root count, and a forced even local factor at two.
In particular, the degree bound on non-even factors follows from
the obligatory block and does not require a separate quintic-factor
argument. The unique-quartic parity restriction follows from the same
normal form and its sign equation.

No global squarefreeness theorem is needed: reflection preserves
multiplicities, the normal form retains them, and the two negative
roots used to bound the even quadratics are simple.
Nor is the complete classification of all local factors required.
The local proof selects one factor with sufficient degree and,
when needed, sufficient control of its constant square class.

Table~\ref{tab:excluded} therefore records $54$ proved impossibilities,
whereas Table~\ref{tab:remaining} records $22$ possibilities not excluded
by the present argument. The latter include patterns with rational
roots, quadratic factors, and cubic factors. Their compatibility with
the stated conditions does not imply their existence for this family.
Excluding them requires additional arithmetic information.
The argument proves neither pointwise irreducibility of the third
polynomial nor nonexistence of perfect cuboids.

\appendix
\section{Exact verification script and transcript}\label{app:certificate}

The proof in the main text is independent of a bounded parameter
search. The following script uses Python and SymPy~\cite{SymPy} with exact
rational arithmetic. The script has been verified with Python~3.12.14
and SymPy~1.14.0. It checks the sextic identities, the shifted
coefficients, the congruence divisibilities used in the local proof,
and the signed normal-form enumeration. The partition list is generated
exhaustively and checked against two independent counting procedures;
all reported counts and lists are computed. Newton polygon deductions,
Hensel lifting, and the field and norm arguments are proved in the text.
The partition enumeration concerns abstract degree data and is finite;
it does not sample the parameters $a,b,u$.

\subsection*{Script: Complete exact verification}
\noindent\textbf{Code.}
\begin{code}[breakanywhere=true]
"""
Exact identities and exhaustive degree-pattern verification.
Requires Python 3.8+ and SymPy. All counts and lists are computed.
"""

from collections import Counter
from dataclasses import asdict, dataclass
from itertools import product
from math import prod
import sympy as s

def zero(expression):
    assert s.cancel(expression) == 0

def integral_coefficients(expression, *variables):
    polynomial = s.Poly(s.expand(expression), *variables, domain=s.QQ)
    assert all(s.denom(coefficient) == 1
               for coefficient in polynomial.coeffs())

x, z, A, B, U, V = s.symbols("x z A B U V")
F = (((x-A)*(x+B)*(x+U))**2 + ((x-B)*(x+A)*(x+U))**2 - ((x-U)*(x+A)*(x+B))**2)
S, Tpar = A+B, A*B
expanded_F = (x**6+(6*U-2*S)*x**5
    +(S*S+4*U*S+U*U-14*Tpar)*x**4
    +(6*U*S*S-20*U*Tpar-2*S*(U*U+Tpar))*x**3
    +(Tpar*Tpar+U*U*S*S-14*Tpar*U*U+4*U*S*Tpar)*x*x
    +(6*U*Tpar*Tpar-2*U*U*S*Tpar)*x+U*U*Tpar*Tpar)
zero(F-expanded_F)
lam = s.symbols("lam")
zero(F.subs({x:lam**2*x, A:lam**2*A,
             B:lam**2*B, U:lam**2*U}, simultaneous=True)-lam**12*F)
zero(s.Poly(F, x).LC()-1)
zero(F.subs(x, 0)-(A*B*U)**2)
zero(F.subs(x, -A)-4*A*A*(A-B)**2*(A-U)**2)
zero(F.subs(x, -B)-4*B*B*(A-B)**2*(B-U)**2)
zero(F.subs(x, -U)+4*U*U*(A-U)**2*(B-U)**2)
RA, RB, RU = [(x-W)/(x+W) for W in (A, B, U)]
zero(F-(x+A)**2*(x+B)**2*(x+U)**2*(RA**2+RB**2-RU**2))

# The real-root identity and its derivative, valid in all three orders.
RVz, RUz = (z-V)/(z+V), (z-U)/(z+U)
qV = RVz/RUz
zero(s.diff(qV, z) - 2*(V-U)*(z*z+U*V)/((z+V)**2*(z-U)**2))
RAz, RBz = (z-A)/(z+A), (z-B)/(z+B)
L = (RUz/RAz)**2+(RUz/RBz)**2
zero(F.subs(x, -z)
     - (z-A)**2*(z-B)**2*(z+U)**2*(L-1))
# The sign of q_V-1 supplies the outside-interval exclusions.
zero(qV-1-2*z*(U-V)/((z+V)*(z-U)))

# Exact expansion near the distinguished point x=-U.
q, C, D = s.symbols("q C D")
FCD = s.expand(F.subs({A:U+C, B:U+D}))
shifted = s.expand(FCD.subs(x, q-U))
expected = [
    -4*C*C*D*D*U*U,
    4*C*D*U*(C*D-2*C*U-2*D*U),
    C*D*(C*D+12*C*U+12*D*U-16*U*U),
    -2*C*D*(C+D-16*U),
    C*C-12*C*D+4*C*U+D*D+4*D*U+4*U*U,
    -2*(C+D+2*U),
    1,
]
zero(shifted-sum(value*q**j for j, value in enumerate(expected)))

# The shared-valuation quartic before the integral perturbation.
E = ((x-U)**2*(x+U+C)**2
     - (x-U-C)**2*(x+U)**2
     - (x+U+C)**2*(x+U)**2)
g = (q**4-2*C*q**3+C*(C+12*U)*q*q
     +4*U*C*(C-2*U)*q-4*U*U*C*C)
zero(E.subs(x, q-U)+g)

# The unique-valuation quartic and its root-field identity.
Q = z*z*(z+2)**2-2*(z*z+2*z+2)**2
zero(Q+z**4+4*z**3+12*z*z+16*z+8)
R = (x-U)/(x+U)
zero((1+R)**2/(R**2-1)+x/U)

# A short exact formula for the outer two-root preparation.
Y = s.symbols("Y")
H = s.cancel(FCD.subs(x, U+4*Y)/(256*U**4))
Hshort = (
    ((Y-C/4)*(1+(D+4*Y)/(2*U))*(1+2*Y/U))**2
    + ((Y-D/4)*(1+(C+4*Y)/(2*U))*(1+2*Y/U))**2
    - (Y*(1+(C+4*Y)/(2*U))*(1+(D+4*Y)/(2*U)))**2
)
zero(H-Hshort)
c, d = s.symbols("c d")
Hunit = s.Poly(s.expand(H.subs({C:8*c, D:8*d})), Y)
# Multiplication by U^4 preserves integrality for an odd unit U.
for j, divisor in enumerate([4, 4, 1, 8, 8, 32, 16]):
    integral_coefficients(Hunit.nth(j)*U**4/divisor, U, c, d)
integral_coefficients((Hunit.nth(2)-1)*U**4/4, U, c, d)
b, e = s.symbols("b e")
outer_factor = (q-2*U)**2+4*b*(q-2*U)+16*e
zero(s.resultant(Y*Y+b*Y+e, q-2*U-4*Y, Y)-outer_factor)
zero(outer_factor.subs(q, 0)-(4*U*U-8*U*b+16*e))

# Equal valuation of C and D: the product identity makes the mod-8
# congruence transparent.  Put C=h*c, D=h*d, with v_2(h)>=3.
h, eta = s.symbols("h eta")
T = s.cancel(FCD.subs({x:-U+h*z, C:h*c, D:h*d})
             /(4*U*U*h**4))
Tshort = (
    (z*(z+d)*(1+h*(c-z)/(2*U)))**2
    + (z*(z+c)*(1+h*(d-z)/(2*U)))**2
    - ((z+c)*(z+d)*(1-h*z/(2*U)))**2
)
zero(T-Tshort)
en, jn = c*d, (c+d)/2
T0 = z**4-en**2-4*en*z*z-4*en*jn*z
zero(Tshort.subs(h, 0)-T0)
# h=8*eta proves every coefficient of T-T0 is divisible by eight
# over Z_2[U,U^{-1},c,d]; this is not a sample parameter check.
integral_coefficients((T-T0).subs(h, 8*eta)*U*U/8,
                      U, eta, c, d, z)

t, p, constant = s.symbols("t p constant")
zero(s.discriminant(t**4+p*t*t+constant, t)
     -16*constant*(p*p-4*constant)**2)
print("PASS: exact polynomial, real-root, and dyadic coefficient identities.")

# Exhaustive inventory and independent completeness checks.
DEGREE = 12
HALF_DEGREE = DEGREE // 2
MAX_EVEN_QUADRATICS = 2
NEGATIVE_ROOTS_OF_SEXTIC = 2

class VerificationError(Exception):
    pass

def require(condition, message):
    # Deliberately not an assert: these checks also run with python -O.
    if not condition:
        raise VerificationError(message)

def partitions(total, minimum=1):
    if total == 0:
        yield ()
        return
    for first in range(minimum, total + 1):
        for tail in partitions(total - first, first):
            yield (first,) + tail

def count_partitions_by_length(total):
    counts = [[0] * (total + 1) for _ in range(total + 1)]
    counts[0][0] = 1
    for n in range(1, total + 1):
        for k in range(1, n + 1):
            counts[n][k] = counts[n - 1][k - 1] + counts[n - k][k]
    return counts[total]

def count_partitions_by_product(total):
    coefficients = [1] + [0] * total
    for part in range(1, total + 1):
        for n in range(part, total + 1):
            coefficients[n] += coefficients[n - part]
    return coefficients[total]

def checked_partitions(total):
    for n in range(total + 1):
        generated = list(partitions(n))
        require(len(generated) == len(set(generated)),
                'Repeated partition at degree {}'.format(n))
        for pattern in generated:
            require(sum(pattern) == n and all(d > 0 for d in pattern),
                    'Invalid partition: {!r}'.format(pattern))
            require(tuple(sorted(pattern)) == pattern,
                    'Unordered partition: {!r}'.format(pattern))
        observed = Counter(map(len, generated))
        counted = count_partitions_by_length(n)
        require(all(observed[k] == counted[k] for k in range(n + 1)),
                'Counts by length disagree at degree {}'.format(n))
        require(len(generated) == sum(counted) == count_partitions_by_product(n),
                'Independent partition counts disagree at degree {}'.format(n))
    # Convert to descending notation only AFTER the completeness checks.
    return tuple(sorted((tuple(reversed(p)) for p in generated), reverse=True))

@dataclass(frozen=True)
class Profile:
    k: int
    ds: tuple
    signs: tuple
    hs: tuple
    negative_roots: tuple

    @property
    def m(self):
        return sum(self.hs)

def reflection_assignments(pattern):
    multiplicities = Counter(pattern)
    degrees = sorted(multiplicities, reverse=True)
    choices = []
    for degree in degrees:
        count = multiplicities[degree]
        if degree % 2:
            if count % 2:
                return
            choices.append((0,))
        else:
            limit = min(count, MAX_EVEN_QUADRATICS) if degree == 2 else count
            choices.append(tuple(f for f in range(limit + 1) if (count - f) % 2 == 0))
    for fixed_counts in product(*choices):
        k, ds, hs = 0, [], []
        for degree, fixed in zip(degrees, fixed_counts):
            if degree == 2:
                k = fixed
            elif fixed:
                ds.extend([degree // 2] * fixed)
            hs.extend([degree] * ((multiplicities[degree] - fixed) // 2))
        yield k, tuple(ds), tuple(hs)

def real_signature(ds, signs, k):
    for negatives in product(*(range(d + 1) for d in ds)):
        if (sum(negatives) == NEGATIVE_ROOTS_OF_SEXTIC - k
                and all(sign == (-1) ** (d + nu)
                        for d, sign, nu in zip(ds, signs, negatives))):
            return tuple(negatives)
    return None

def normal_form_profiles(pattern):
    profiles = []
    for k, ds, hs in reflection_assignments(pattern):
        m = sum(hs)
        require(k + sum(ds) + m == HALF_DEGREE,
                'Degree bookkeeping failed for {!r}'.format(pattern))
        if not ds:
            continue
        for signs in product((-1, 1), repeat=len(ds)):
            if prod(signs) != (-1) ** m:
                continue
            obligatory_factor = any(d >= 3 or (d == 2 and sign == 1)
                                    for d, sign in zip(ds, signs))
            if not obligatory_factor:
                continue
            negatives = real_signature(ds, signs, k)
            # This independently verifies that the exact real root count
            # does not further reduce the signed profiles in this inventory.
            require(negatives is not None,
                    'A signed profile fails the real-root test: {!r}'.format(
                        (pattern, k, ds, signs, hs)))
            profiles.append(Profile(k, ds, tuple(signs), hs, negatives))
    return tuple(profiles)

def exclusion_rule(pattern):
    counts = Counter(pattern)
    unpaired = sorted((d for d, count in counts.items() if d % 2 and count % 2),
                      reverse=True)
    if unpaired:
        return 'I', 'odd multiplicity of odd degree(s): ' + ', '.join(map(str, unpaired))
    if not any(d >= 4 and d % 2 == 0 for d in pattern):
        return 'II', 'no even degree >= 4'
    if max(pattern) <= 4 and counts[4] == 1 and counts[2] % 2:
        return 'III', 'unique quartic; odd number of quadratics'
    return None, None

def verify_profile(pattern, profile):
    require(len(profile.ds) == len(profile.signs) == len(profile.negative_roots),
            'Witness vector lengths disagree')
    require(all(sign in (-1, 1) for sign in profile.signs),
            'Witness has an invalid constant-term sign')
    require(all(h > 0 for h in profile.hs), 'Witness has a nonpositive H degree')
    reconstructed = ([2] * profile.k + [2 * d for d in profile.ds]
                     + [h for h in profile.hs for _ in range(2)])
    require(tuple(sorted(reconstructed, reverse=True)) == pattern,
            'Witness does not reconstruct {}'.format(pattern))
    require(0 <= profile.k <= MAX_EVEN_QUADRATICS and bool(profile.ds),
            'Invalid quadratic bound or missing obligatory factor')
    require(all(d >= 2 for d in profile.ds), 'Invalid higher even degree')
    require(profile.k + sum(profile.ds) + profile.m == HALF_DEGREE,
            'Witness degree equation failed')
    require(prod(profile.signs) == (-1) ** profile.m, 'Witness sign equation failed')
    require(any(d >= 3 or (d == 2 and sign == 1)
                for d, sign in zip(profile.ds, profile.signs)),
            'Witness obligatory-factor condition failed')
    require(sum(profile.negative_roots) + profile.k == NEGATIVE_ROOTS_OF_SEXTIC,
            'Witness real-root count failed')
    require(all(0 <= nu <= d and sign == (-1) ** (d + nu)
                for d, sign, nu in zip(profile.ds, profile.signs, profile.negative_roots)),
            'Witness real-root signs failed')

def audit():
    universe = checked_partitions(DEGREE)
    records = []
    for pattern in universe:
        profiles = normal_form_profiles(pattern)
        stage, reason = exclusion_rule(pattern)
        require(bool(profiles) == (stage is None),
                'The two classification methods disagree for {!r}'.format(pattern))
        for profile in profiles:
            verify_profile(pattern, profile)
        records.append({
            'pattern': pattern,
            'status': ('irreducible' if len(pattern) == 1
                       else 'retained' if profiles else 'excluded'),
            'exclusion_stage': stage,
            'reason': reason,
            'witness': asdict(profiles[0]) if profiles else None,
        })
    trivial = [r for r in records if r['status'] == 'irreducible']
    retained = [r for r in records if r['status'] == 'retained']
    excluded = [r for r in records if r['status'] == 'excluded']
    require(len(trivial) == 1 and trivial[0]['pattern'] == (DEGREE,),
            'The one-part pattern has not been handled correctly')
    require(len(records) == len(trivial) + len(retained) + len(excluded),
            'The final classes do not cover the universe')
    stages = Counter(r['exclusion_stage'] for r in excluded)
    summary = {
        'degree': DEGREE,
        'all_partitions': len(universe),
        'partition_count_by_length': dict(sorted(Counter(map(len, universe)).items())),
        'irreducible_patterns': len(trivial),
        'nontrivial_patterns': len(universe) - len(trivial),
        'excluded_patterns': len(excluded),
        'retained_patterns': len(retained),
        'excluded_by_stage': {stage: stages[stage] for stage in ('I', 'II', 'III')},
        'retained_by_largest_degree': dict(sorted(
            Counter(r['pattern'][0] for r in retained).items(), reverse=True)),
    }
    return {'summary': summary, 'records': records}

def format_pattern(pattern):
    return '+'.join(map(str, pattern))

result = audit()
summary, records = result["summary"], result["records"]
print("Partition counts by number of factors:")
print(", ".join("({}: {})".format(length, count)
      for length, count in summary["partition_count_by_length"].items()))
print("All partitions: {}".format(summary["all_partitions"]))
print("Nontrivial partitions: {} - {} = {}".format(
    summary["all_partitions"], summary["irreducible_patterns"],
    summary["nontrivial_patterns"]))
stages = summary["excluded_by_stage"]
print("Excluded by I, II, III: {} + {} + {} = {}".format(
    stages["I"], stages["II"], stages["III"],
    summary["excluded_patterns"]))
print("Nontrivial = excluded + retained: {} = {} + {}".format(
    summary["nontrivial_patterns"], summary["excluded_patterns"],
    summary["retained_patterns"]))
print("Retained counts by largest degree:")
print(", ".join("({}: {})".format(degree, count)
      for degree, count in summary["retained_by_largest_degree"].items()))
print("Complete degree inventory:")
for index, record in enumerate(records, 1):
    status = ("ONE-PART" if record["status"] == "irreducible"
              else record["status"].upper())
    if record["exclusion_stage"]:
        status += " ({}): {}".format(
            record["exclusion_stage"], record["reason"])
    print("{:2d}. {:25s} {}".format(
        index, format_pattern(record["pattern"]), status))
print("Formal compatibility witnesses (k, ds, eps, hs, nu):")
for record in records:
    if record["status"] != "retained":
        continue
    w = record["witness"]
    print("{}: {}".format(format_pattern(record["pattern"]),
          (w["k"], w["ds"], w["signs"],
           w["hs"], w["negative_roots"])))
\end{code}

\noindent\textbf{Transcript.}
\begin{term}
PASS: exact polynomial, real-root, and dyadic coefficient identities.
Partition counts by number of factors:
(1: 1), (2: 6), (3: 12), (4: 15), (5: 13), (6: 11), (7: 7), (8: 5), (9: 3), (10: 2), (11: 1), (12: 1)
All partitions: 77
Nontrivial partitions: 77 - 1 = 76
Excluded by I, II, III: 37 + 14 + 3 = 54
Nontrivial = excluded + retained: 76 = 54 + 22
Retained counts by largest degree:
(10: 2), (8: 4), (6: 8), (4: 8)
Complete degree inventory:
 1. 12                        ONE-PART
 2. 11+1                      EXCLUDED (I): odd multiplicity of odd degree(s): 11, 1
 3. 10+2                      RETAINED
 4. 10+1+1                    RETAINED
 5. 9+3                       EXCLUDED (I): odd multiplicity of odd degree(s): 9, 3
 6. 9+2+1                     EXCLUDED (I): odd multiplicity of odd degree(s): 9, 1
 7. 9+1+1+1                   EXCLUDED (I): odd multiplicity of odd degree(s): 9, 1
 8. 8+4                       RETAINED
 9. 8+3+1                     EXCLUDED (I): odd multiplicity of odd degree(s): 3, 1
10. 8+2+2                     RETAINED
11. 8+2+1+1                   RETAINED
12. 8+1+1+1+1                 RETAINED
13. 7+5                       EXCLUDED (I): odd multiplicity of odd degree(s): 7, 5
14. 7+4+1                     EXCLUDED (I): odd multiplicity of odd degree(s): 7, 1
15. 7+3+2                     EXCLUDED (I): odd multiplicity of odd degree(s): 7, 3
16. 7+3+1+1                   EXCLUDED (I): odd multiplicity of odd degree(s): 7, 3
17. 7+2+2+1                   EXCLUDED (I): odd multiplicity of odd degree(s): 7, 1
18. 7+2+1+1+1                 EXCLUDED (I): odd multiplicity of odd degree(s): 7, 1
19. 7+1+1+1+1+1               EXCLUDED (I): odd multiplicity of odd degree(s): 7, 1
20. 6+6                       RETAINED
21. 6+5+1                     EXCLUDED (I): odd multiplicity of odd degree(s): 5, 1
22. 6+4+2                     RETAINED
23. 6+4+1+1                   RETAINED
24. 6+3+3                     RETAINED
25. 6+3+2+1                   EXCLUDED (I): odd multiplicity of odd degree(s): 3, 1
26. 6+3+1+1+1                 EXCLUDED (I): odd multiplicity of odd degree(s): 3, 1
27. 6+2+2+2                   RETAINED
28. 6+2+2+1+1                 RETAINED
29. 6+2+1+1+1+1               RETAINED
30. 6+1+1+1+1+1+1             RETAINED
31. 5+5+2                     EXCLUDED (II): no even degree >= 4
32. 5+5+1+1                   EXCLUDED (II): no even degree >= 4
33. 5+4+3                     EXCLUDED (I): odd multiplicity of odd degree(s): 5, 3
34. 5+4+2+1                   EXCLUDED (I): odd multiplicity of odd degree(s): 5, 1
35. 5+4+1+1+1                 EXCLUDED (I): odd multiplicity of odd degree(s): 5, 1
36. 5+3+3+1                   EXCLUDED (I): odd multiplicity of odd degree(s): 5, 1
37. 5+3+2+2                   EXCLUDED (I): odd multiplicity of odd degree(s): 5, 3
38. 5+3+2+1+1                 EXCLUDED (I): odd multiplicity of odd degree(s): 5, 3
39. 5+3+1+1+1+1               EXCLUDED (I): odd multiplicity of odd degree(s): 5, 3
40. 5+2+2+2+1                 EXCLUDED (I): odd multiplicity of odd degree(s): 5, 1
41. 5+2+2+1+1+1               EXCLUDED (I): odd multiplicity of odd degree(s): 5, 1
42. 5+2+1+1+1+1+1             EXCLUDED (I): odd multiplicity of odd degree(s): 5, 1
43. 5+1+1+1+1+1+1+1           EXCLUDED (I): odd multiplicity of odd degree(s): 5, 1
44. 4+4+4                     RETAINED
45. 4+4+3+1                   EXCLUDED (I): odd multiplicity of odd degree(s): 3, 1
46. 4+4+2+2                   RETAINED
47. 4+4+2+1+1                 RETAINED
48. 4+4+1+1+1+1               RETAINED
49. 4+3+3+2                   EXCLUDED (III): unique quartic; odd number of quadratics
50. 4+3+3+1+1                 RETAINED
51. 4+3+2+2+1                 EXCLUDED (I): odd multiplicity of odd degree(s): 3, 1
52. 4+3+2+1+1+1               EXCLUDED (I): odd multiplicity of odd degree(s): 3, 1
53. 4+3+1+1+1+1+1             EXCLUDED (I): odd multiplicity of odd degree(s): 3, 1
54. 4+2+2+2+2                 RETAINED
55. 4+2+2+2+1+1               EXCLUDED (III): unique quartic; odd number of quadratics
56. 4+2+2+1+1+1+1             RETAINED
57. 4+2+1+1+1+1+1+1           EXCLUDED (III): unique quartic; odd number of quadratics
58. 4+1+1+1+1+1+1+1+1         RETAINED
59. 3+3+3+3                   EXCLUDED (II): no even degree >= 4
60. 3+3+3+2+1                 EXCLUDED (I): odd multiplicity of odd degree(s): 3, 1
61. 3+3+3+1+1+1               EXCLUDED (I): odd multiplicity of odd degree(s): 3, 1
62. 3+3+2+2+2                 EXCLUDED (II): no even degree >= 4
63. 3+3+2+2+1+1               EXCLUDED (II): no even degree >= 4
64. 3+3+2+1+1+1+1             EXCLUDED (II): no even degree >= 4
65. 3+3+1+1+1+1+1+1           EXCLUDED (II): no even degree >= 4
66. 3+2+2+2+2+1               EXCLUDED (I): odd multiplicity of odd degree(s): 3, 1
67. 3+2+2+2+1+1+1             EXCLUDED (I): odd multiplicity of odd degree(s): 3, 1
68. 3+2+2+1+1+1+1+1           EXCLUDED (I): odd multiplicity of odd degree(s): 3, 1
69. 3+2+1+1+1+1+1+1+1         EXCLUDED (I): odd multiplicity of odd degree(s): 3, 1
70. 3+1+1+1+1+1+1+1+1+1       EXCLUDED (I): odd multiplicity of odd degree(s): 3, 1
71. 2+2+2+2+2+2               EXCLUDED (II): no even degree >= 4
72. 2+2+2+2+2+1+1             EXCLUDED (II): no even degree >= 4
73. 2+2+2+2+1+1+1+1           EXCLUDED (II): no even degree >= 4
74. 2+2+2+1+1+1+1+1+1         EXCLUDED (II): no even degree >= 4
75. 2+2+1+1+1+1+1+1+1+1       EXCLUDED (II): no even degree >= 4
76. 2+1+1+1+1+1+1+1+1+1+1     EXCLUDED (II): no even degree >= 4
77. 1+1+1+1+1+1+1+1+1+1+1+1   EXCLUDED (II): no even degree >= 4
Formal compatibility witnesses (k, ds, eps, hs, nu):
10+2: (1, (5,), (1,), (), (1,))
10+1+1: (0, (5,), (-1,), (1,), (2,))
8+4: (0, (4, 2), (-1, -1), (), (1, 1))
8+2+2: (0, (4,), (1,), (2,), (2,))
8+2+1+1: (1, (4,), (-1,), (1,), (1,))
8+1+1+1+1: (0, (4,), (1,), (1, 1), (2,))
6+6: (0, (3, 3), (-1, -1), (), (0, 2))
6+4+2: (1, (3, 2), (-1, -1), (), (0, 1))
6+4+1+1: (0, (3, 2), (-1, 1), (1,), (0, 2))
6+3+3: (0, (3,), (-1,), (3,), (2,))
6+2+2+2: (1, (3,), (1,), (2,), (1,))
6+2+2+1+1: (0, (3,), (-1,), (2, 1), (2,))
6+2+1+1+1+1: (1, (3,), (1,), (1, 1), (1,))
6+1+1+1+1+1+1: (0, (3,), (-1,), (1, 1, 1), (2,))
4+4+4: (0, (2,), (1,), (4,), (2,))
4+4+2+2: (0, (2, 2), (1, 1), (2,), (0, 2))
4+4+2+1+1: (1, (2, 2), (-1, 1), (1,), (1, 0))
4+4+1+1+1+1: (0, (2, 2), (1, 1), (1, 1), (0, 2))
4+3+3+1+1: (0, (2,), (1,), (3, 1), (2,))
4+2+2+2+2: (0, (2,), (1,), (2, 2), (2,))
4+2+2+1+1+1+1: (0, (2,), (1,), (2, 1, 1), (2,))
4+1+1+1+1+1+1+1+1: (0, (2,), (1,), (1, 1, 1, 1), (2,))
\end{term}

\section*{AI Assistance}
OpenAI GPT-6.1 Sol assisted with the development of the argument, symbolic
calculations, verification code, and the preparation of the exposition.
The author is responsible for the mathematical content and the final
presentation.


\begingroup
\footnotesize
\begin{thebibliography}{99}

\bibitem{SharipovCuboids}
R.~A.~Sharipov,
\newblock \emph{Perfect cuboids and irreducible polynomials},
\newblock Ufa Math.\ J.\ \textbf{4} (2012), no.~1, 153--160.

\bibitem{SharipovThird}
R.~A.~Sharipov,
\newblock \emph{A note on the third cuboid conjecture. Part~I},
\newblock arXiv:1203.2567 [math.NT] (2012).

\bibitem{AndrewsPartitions}
G.~E.~Andrews,
\newblock \emph{\href{https://doi.org/10.1017/CBO9780511608650}{The Theory of Partitions}},
\newblock Cambridge Mathematical Library,
\newblock Cambridge University Press, Cambridge, 1998.

\bibitem{Serre}
J.-P.~Serre,
\newblock \emph{\href{https://doi.org/10.1007/978-1-4757-5673-9}{Local Fields}},
\newblock Graduate Texts in Mathematics, vol.~67,
\newblock Springer-Verlag, New York--Berlin, 1979.

\bibitem{Robert}
A.~M.~Robert,
\newblock \emph{\href{https://doi.org/10.1007/978-1-4757-3254-2}{A Course in $p$-adic Analysis}},
\newblock Graduate Texts in Mathematics, vol.~198,
\newblock Springer-Verlag, New York, 2000.

\bibitem{Berger}
L.~Berger,
\newblock \emph{\href{https://doi.org/10.17879/59019521586}{The Weierstrass preparation theorem and resultants of $p$-adic power series}},
\newblock M\"unster J.\ Math.\ \textbf{14} (2021), 155--163.

\bibitem{SymPy}
A.~Meurer, C.~P.~Smith, M.~Paprocki, et~al.,
\newblock \emph{\href{https://doi.org/10.7717/peerj-cs.103}{SymPy: symbolic computing in Python}},
\newblock PeerJ Comput.\ Sci.\ \textbf{3} (2017), e103.

\end{thebibliography}
\endgroup
\end{document}